\documentclass[10pt,fleqn,a4paper,oneside,draft]{article}
\usepackage[english]{style}

\usepackage[all,ps,tips,dvips]{xy}
\SelectTips{eu}{10}

\def\dash---{\thinspace---\hskip.16667em\relax}

\let\to\rightarrow

\let\l\lambda
\let\r\rho
\let\*\bullet
\let\ot\otimes
\let\cot\square
\def\place{\mathord{?}}

\def\sub#1{{}_{#1}}

\let\category\mathsf
\def\C{\category{C}}
\def\D{\category{D}}
\def\Mod{\category{Mod}}
\def\M{\mathcal{M}}

\def\HH{H\mskip-2mu H}

\def\Ae{A\!^e}

\def\.{\mskip1mu}

\def\HB{{}_H^H\Mod{}_H^H}

\def\ZZ{\mathbb{Z}}

\makeatletter
\def\operator@font{\sf}
\makeatother

\DeclareMathOperator{\End}{end} 
\DeclareMathOperator{\Ext}{Ext}
\DeclareMathOperator{\obj}{obj}

\def\xycenter#1{\xy*[c]\xybox{#1}\endxy}

\def\xytriangle#1#2#3#4#5#6{
  \begin{xy}
  +*[c]{\object{
    \begin{xy}
      <1.5cm,0cm>:
      *+[l]{#1}="a",
      +a(-60) *+{#6}="c",
      +a(60) *+[r]{#4}="b"
      \POS"a" \ar"b" ^-{#2}
      \POS"a" \ar"c" _-{#3}
      \POS"b" \ar"c" ^-{#5}
    \end{xy}
    }}
  \end{xy}
}

\makeatletter
\def\@biblabel#1{#1.}
\makeatother

\title{The Hilton-Eckmann argument for cup products}
\author{%
Mariano Suarez-Alvarez\thanks{Departamento de Matem\'atica,
Facultad de Ciencias Exactas y Naturales, Universidad
de Buenos Aires. Ciudad Universitaria. Pabellón~I, Buenos Aires (1428) 
Argentina. e-mail: \texttt{mariano@dm.uba.ar}.\newline
\hspace*{\parindent} 
This work was supported by a grant from \textsc{UBACyT} X062, the international
cooperation project \textsc{SECyT-ECOS} A98E05, and a \textsc{CoNICET}
scholarship.}
}
\date{}

\begin{document}

\maketitle

\begin{abstract}
We present a simple extension of the classical Hilton-Eckmann argument which
proves that the endomorphism monoid of the unit object in a monoidal category is
commutative. It allows us to recover in a uniform way well-known results on the
graded-commutativity of cup products defined on the cohomology theories attached
to various algebraic structures, as well as some more recent results.
\end{abstract}

\section*{Introduction}

In their paper~\cite{EckmannHilton} on group objects in general categories,
B.~Eckmann and P.~J.~Hilton show that two monoid structures on a set which are
such that one of them is a homomorphism for the other actually coincide and
are commutative. 

Interestingly, this situation does occur in nature: if $X$ and $Y$ are
respectively a co-$H$-space and an $H$-space, then the set of based homotopy
classes of maps $[X,Y]$ is a monoid in two ways, and naturality implies that
each of the corresponding products on $[X,Y]$ is a homomorphism for the other.
The result just quoted implies then that these two structures are equal, and
that $[X,Y]$ is commutative. Instances of this are the well-known facts that
both the fundamental group $\pi_1(G,e)$ of a topological group $G$ and the
higher homotopy groups $\pi_p(X,x_0)$ for $p>1$ of any space $X$, are abelian.

Another situation in which the same argument applies is the following. Let $\C$
be a monoidal category with unit object $e\in\obj\C$. The set $\End_\C(e)$ of
endomorphisms of $e$ in $\C$ is of course an associative monoid with respect to
composition. It turns out that it is always commutative. Indeed, we can define
on $\End_\C(e)$ a convolution product and verify that it is a homomorphism with
respect to composition. It follows that composition and convolution in
$\End_\C(e)$ are equal and commutative. We refer to Ch.~Kassel's
book~\cite{Kassel}, section \textsc{XI}.2, for details.

\medskip

The purpose of this note is to present a simple extension of the argument of
Eckmann and Hilton which can be used to show that products defined on various
cohomology theories are commutative. It applies to the cohomology of groups, to
the Hochschild cohomology of associative algebras, to the Cartier cohomology of
coalgebras, and to other somewhat more exotic theories, as the cohomology theory
introduced by M.~Gerstenhaber and S.~D.~Schack in~\cite{GerstenhaberSchack} for Hopf bimodules
over a Hopf algebra.

\medskip

In the next section we recall the definition of a monoidal category so as to fix
the notation, we extend it in a natural way in order to adapt it to a
``derived'' context, and state and prove our theorem. In section 2 we indicate
how our result allows us to establish easily the graded-commutativity of the
product on the usual cohomology theories.

\section{Definitions and the theorem}

\paragraph In what follows all categories and functors are implicitly assumed to
be additive.

\paragraph A \emph{monoidal category} is a $6$-tuple $(\C,\ot,e,a,l,r)$ in which
$\C$ is a category, $?\ot?:\C\times\C\to\C$ a bifunctor, $e\in\obj\C$ an object,
$a:(\place\ot\place)\ot\place\to\place\ot(\place\ot\place)$ an isomorphism of
functors $\C\times\C\times\C\to\C$, and in which $l:e\ot?\to?$ and $r:?\ot
e\to?$ are isomorphisms of functors $\C\to\C$, which are such that for each
choice of objects $x$,~$y$,~$z$ and~$w\in\obj\C$, the following diagrams
commute:
  \[
  \xycenter{\xymatrix{
    (x\ot(y\ot z))\ot w \ar[dd]_-a
      & ((x\ot y)\ot z)\ot w) \ar[l]_-{a\ot1} \ar[d]^-{a}
      \\
    {}
      & (x\ot y)\ot(z\ot w) \ar[d]^-{a}
      \\
    x\ot((y\ot z)\ot w) \ar[r]^{1\ot a}
      & x\ot(y\ot(z\ot w))
  }}
  \quad
  \xytriangle{(x\ot e)\ot y}{a}{r\ot1}{x\ot(e\ot y)}{1\ot l}{x\ot y}
  \]

\paragraph It follows easily from the definition that $l=r:e\ot e\to e$; see for
example \cite{Kassel}, lemma~\textsc{XI}.2.3. The proof of this fact involves the
associativity constraint $a$, and it is interesting to notice that this is the
only r\^ole played by $a$ in what follows.

\paragraph\label{def:smc} A \emph{suspended monoidal category} is a $9$-tuple
$(\C,\ot,e,a,l,r,T,\l,\r)$ such that $(\C,\ot,e,a,r,l)$ is a monoidal category,
$T:\C\to\C$ is an automorphism, and $\l:\place\ot T\place\to T(\place\ot\place)$
and $\r:T\place\ot\place\to T(\place\ot\place)$ are isomorphisms of functors
$\C\times\C\to\C$ such that for each pair of objects $x$ and~$y\in\obj\C$, the
following two diagrams commute
  \[
  \xymatrix{
    e\ot Tx \ar[r]^-l \ar[d]_-\l
      & Tx \ar[d]^1
      \\
    T(e\ot x) \ar[r]^-{T\.l}
      & Tx
    } 
  \qquad
  \xymatrix{
    Tx\ot e \ar[r]^-{r} \ar[d]_-\r
      & Tx \ar[d]^1
      \\
    T(x\ot e) \ar[r]^-{Tr}
      & Tx
    }
  \]
and the following diagram anti-commutes
  \[
  \xymatrix{
    Tx\ot Ty \ar[r]^-\r \ar[d]_-\l \ar@{}[rd]|-{\mbox{$-1$}}
      & T(x\ot Ty) \ar[d]^-{T\.\l}
      \\
    T(Tx\ot y) \ar[r]^-{T\rho}
      & T^2(x\ot y)
    }
  \]

\paragraph Given a suspended monoidal category as in~\pref{def:smc}, we put
$\l_0=\r_0=1:x\ot y\to x\ot y$, and for each $p>0$, 
  \begin{gather*}
  \l_p = T^{p-1}\l\circ T^{p-2}\l\circ\cdots\circ\l : x\ot T^py \to T^p(x\ot y), \\
  \r_p = T^{p-1}\r\circ T^{p-2}\r\circ\cdots\circ\r : T^px\ot y \to T^p(x\ot y), \\
\intertext{and}
  \l_{-p} = T^{-p}\l_p^{-1}:x\ot T^{-p}y\to T^p(x\ot y), \\
  \r_{-p} = T^{-p}\r_p^{-1}:T^{-p}x\ot y\to T^p(x\ot y). 
  \end{gather*}
Observe that this last two equations do make sense: for example, since
$\l_p:x\ot T^pT^{-p}y\to T^p(x\ot T^{-p}y)$, we have $\l^{-1}_p:T^p(x\ot
T^{-p}y)\to x\ot T^pT^{-p}y$ so $T^{-p}\l^{-1}_p:x\ot T^{-p}y\to T^{-p}(x\ot y)$.

\paragraph We obtain in this way isomorphisms of functors $\l_p:\place\ot
T^p\place\to T^p(\place\ot\place)$ and $\r_p:T^p\place\ot\place\to
T^p(\place\ot\place)$ for all $p\in\ZZ$ such that we have
  \begin{gather*}
  l = T^pl\circ\l_p:e\ot T^p\place\to T^p\place \\
  r = T^pr\circ\r_p:T^p\place\ot e\to T^p\place \\
  T^p\l_q \circ \r_p = (-1)^{pq}\;T^q\r_p\circ\l_q
                :T^p\place\ot T^q\place\to T^{p+q}(\place\ot\place)
  \end{gather*}
for all $p$,~$q\in\ZZ$.

\begin{Theorem}\label{thm}
Let $(\C,\ot,e,a,l,r,T,\l,\r)$ be a suspended monoidal category, and put 
  \[
  \End^T_C(e)=\bigoplus_{p\in\ZZ}\hom_C(e,T^pe). 
  \]
If $f:e\to T^pe$ and $g:e\to T^qe$, define $f\cdot g = T^qf\circ g:e\to T^{p+q}e$. 
Then $(\End^T_C(e),\cdot)$ is a commutative graded ring.
\end{Theorem}

\begin{Proof}
The facts that the operation $\cdot$ is associative, that it admits
$1\in\hom_C(e,e)$ as a unit element, and that it is distributive with respect
to addition follow immediately from the corresponding facts about the
composition of morphisms in $\C$. We need only prove then that it is
commutative.

If $f:e\to T^pe$ and $g:e\to T^qe$, we let $f\star g:e\to T^{p+q}e$ be the
composition

  \[\xymatrix@1{
    e \ar[r]^-{r^{-1}}
      & e\ot e \ar[r]^-{f\ot g}
      & T^pe\ot T^qe \ar[r]^-{\r_p}
      & T^p(e\ot T^qe) \ar[r]^-{T^p\l_q}
      & T^{p+q}(e\ot e) \ar[r]^-{T^{p+q}l\;}
      & T^{p+q}e
    }
  \]

Each bounded face in the following diagram commutes:
  \[\xymatrix@C+3pt@R+3pt{
    e \ar[d]_-{f}
      & e\ot e \ar[l]_-{r} \ar[d]_-{f\ot 1} \ar`r[rd]^-{f\ot g}[rd]
      \\
    T^pe \ar[d]_-1
      & T^pe\ot e \ar[l]_-{r} \ar[r]^-{1\ot g} \ar[d]_-{\r_p}
      & T^pe\ot T^qe \ar[d]_-{\r_p}
      \\
    T^pe \ar`d[rd]_-1[rd]
      & T^p(e\ot e) \ar[l]_-{T^pr} \ar[r]^-{T^p(1\ot g)} \ar[d]_-{T^pl}
      & T^p(e\ot T^qe) \ar[d]_-{T^pl} \ar`r[rd]^-{T^p\l_q}[rd]
      \\
    {}
      & T^pe \ar[r]^-{T^pg}
      & T^{p+q}e
      & T^{p+q}(e\ot e)\ar[l]_-{T^{p+q}l}
    }
  \]
so the outer one does, too, and we see that $f\star g=g\cdot f$. On the other
hand, from the $(-1)^{pq}$-commutativity of
  \[\xymatrix@C+3pt@R+3pt{
    e 
                \ar[d]_-{g} 
      & e\ot e  
                \ar@<-3pt>[l]_-{r}
                \ar@<3pt>[l]^-{l} 
                \ar`r[rd]^-{f\ot g}[rd] 
                \ar[d]_-{1\ot g}
      & {} \\
    T^qe 
                \ar[d]_-1
      & e\ot T^qe 
                \ar[l]_-{l} 
                \ar[r]^-{f\ot 1} 
                \ar[d]_-{\l_q}
      & T^pe\ot T^qe 
                \ar[r]^-{\r_p} 
                \ar[d]_-{\l_q} 
                \ar@{}[rd]|-{\quad(-1)^{pq}}
      & T^p(e\ot T^qe) 
                \ar[d]^-{T^p\l_q} 
      & {} \\
    T^qe 
                \ar`d[dr]_-1[dr]
      & T^q(e\ot e) 
                \ar[l]_-{T^ql} 
                \ar[d]_-{T^qr} 
                \ar[r]^-{T^q(f\ot 1)}
      & T^q(T^pe\ot e) 
                \ar[r]^-{T^q\r_p} 
                \ar[d]_-{T^qr}
      & T^{p+q}(e\ot e) 
                \ar[ld]^-{T^{p+q}r} 
                \ar[d]^-{T^{p+q}l}
      & {} \\
    {}
      & T^qe    
                \ar[r]_-{T^qf}
      & T^{p+q}e 
                \ar[r]_-1
      & T^{p+q}e
      & {}
    }
  \]
we see that we also have $f\star g=(-1)^{pq}f\cdot g$. It follows from this, of
course, that~$\cdot$~is a graded-commutative operation.~\qed
\end{Proof}

\section{Applications}

\paragraph\label{p:app} Suppose $(\C,\ot,e,a,l,r)$ is a monoidal category such
that the underlying category $\C$ is an exact category, and such that the tensor
product is an exact functor. Then $\ot$ can be extended naturally to the bounded
derived category $\D^b(\C)$, and it is not difficult to see that $\D^b(\C)$
becomes in this way a suspended monoidal category, whose unit object is simply
$e$ considered as an object of $\D^b(\C)$ in the usual way, and whose suspension
functor is the translation. 

Our theorem implies then that
$\Ext^\*_\C(e,e)=\bigoplus_{p\in\ZZ}\hom_{\D^b(\C)}(e,T^pe)$ is commutative for
the composition product considered in~\pref{thm}, which of course is simply the
Yoneda product.

\paragraph This applies in particular to the abelian monoidal category $\C=\sub
H\Mod$  of left $H$-modules over a Hopf algebra $H$ defined over a field $k$,
with tensor product induced by the tensor product $\ot_k$ of vector spaces
endowed with diagonal action. Since the unit object is $k$ with trivial action,
we conclude that $\Ext_H^\*(k,k)$ is commutative for the Yoneda product.

Dually, the Yoneda algebra $\Ext_H^\*(k,k)$ of self-extensions of $k$ in the
monoidal category ${}^H\Mod$ of left $H$-comodules is commutative.

\paragraph A well-known instance of this situation arises when $H=kG$ is the
group algebra of a group $G$ over a commutative base ring $k$. We recover in
this way the fact that $H^\*(G)=\Ext_{kG}^\*(k,k)$, the group cohomology of $G$
over $k$, is a commutative ring for the Yoneda product. Since this coincides
with the cup product on $H^\*(G)$, we recover the fact that the cup product is
commutative.

\paragraph We remark that all that is really needed in order to be able to
conclude in this way that $\Ext_H^\*(k,k)$ is commutative is a monoidal
structure on $\sub H\Mod$ for which $k$ is the unit object, and this monoidal
structure need not arise from a bialgebra structure on $H$. 

Thus, if $H$ is only a quasi-bialgebra in the sense of
Drinfel'd\dash---see for example~\cite{Kassel}, Chapter~XV---we reach the same
conclusion.

\paragraph The theorem can also be applied in situations in which the tensor
product $\ot$ in the initial datum is not exact but can be extended by
derivation to the derived category. 

For example, if $A$ is a possibly non-commutative $k$-algebra, then the category
$\sub{\Ae}\Mod$ of $A$-bimodules can be endowed with the structure of a monoidal
category with product given by the usual tensor product $\ot_A$ of $A$-modules,
which is in general not an exact functor. It is right exact, though, so it does
admit the left derived functor $\ot_A^L$ as an extension on the derived category
$\D^{-}(\sub{\Ae}\Mod)$ of bounded above complexes, which becomes naturally in
this way a suspended monoidal category\dash---indeed, definition~\pref{def:smc}
is designed to capture precisely this situation. The unit object is clearly
$A\in\obj D^-(\sub{\Ae}\Mod)$. 

Our theorem implies then that the Hochschild cohomology of $A$,
$\HH^\*(A)=\Ext^\*_{\Ae}(A,A)$, is commutative for the Yoneda product. Now,
since the Yoneda product on $\HH^\*(A)$ is the same as the cup product defined
by M.~Gerstenhaber in~\cite{Gerstenhaber1,Gerstenhaber2}, we recover the fact
that this last product is commutative.

We note that this argument can be dualised to obtain a proof that the Yoneda
product on the Cartier cohomology $\HH^\*(C)$ of a coalgebra $C$, or
equivalently, the algebra $\Ext_{C^e}^\*(C,C)$ of self-extensions of $C$ in the
category ${}^C\Mod{}^C$ of $C$-bicomodules is commutative.

\paragraph For our final example, let $H$ be a Hopf algebra over a field $k$,
and let us consider the abelian category $\M=\HB$ of Hopf $H$-bimodules. It can
be made into a monoidal category with product given by the tensor product
$\ot_H$ of $H$-bimodules in such a way that the forgetful functor
$\M\to{}_H\Mod{}_H$ is monoidal. If $M$,~$N\in\M$, then $M\ot_HN$ has its left
module structure induced by that on $M$, its right module structure induced by
that on $N$, and codiagonal left and right comodule structures. We refer to
P.~Schauenburg's paper~\cite{Schauenburg} for details.

The unit object is $H$, and the product is exact because one-sided Hopf
bimodules are free. We are thus in the situation of~\pref{p:app},~and we
conclude that the Yoneda product on the extension algebra $\Ext^\*_\M(H,H)$ is
commutative. 

This result has been obtained in a different way by R.~Taillefer
in~\cite{Taillefer1}, section~4. Since she has shown in~\cite{Taillefer2} that
the Yoneda algebra $\Ext^\*_\M(H,H)$ is isomorphic to the cohomology groups
$H_{GS}^\*(H,H)$ introduced by Gerstenhaber and Schack when these are endowed
with a certain (rather complex) cup product, we see that this last product is
commutative.

We note in passing that the category $\M$ can be endowed with another monoidal
structure, with product induced by the cotensor product $\cot^H$ of
$H$-bicomodules, in such a way that now the forgetful functor
$\M\to{}^H\Mod{}^H$ is monoidal. We do not really obtain anything new, on one
hand because the identity functor $\M\to\M$ is (in a non-trivial way) a monoidal
equivalence $(\M,\ot_H)\cong(\M,\cot^H)$, as shown in~\cite{Schauenburg},
Corollary~6.1, and because the conclusion in~\pref{thm} is actually independent
of the particular monoidal structure under consideration. However, in view of
the equality $f\star g=g\cdot f$ obtained in the course of proof of the theorem,
this shows that the cup product on $\Ext^\*_\M(H,H)$ can be computed using
either the $\star$-product defined in terms of $\ot_H$ or the one defined in
terms of $\cot^H$, and this observation might lead to simplifications when doing
actual computation.


\end{document}